\documentclass[onefignum,onetabnum]{siamart190516}

\usepackage{epsfig}
\usepackage{amsmath}
\usepackage{amssymb}
\usepackage{amsfonts}
\usepackage{graphicx}
\usepackage{epstopdf}
\usepackage{color,soul}
\usepackage{tikz}
\usetikzlibrary{trees}
\usepackage{cleveref}

\newcommand{\order}[1]{{\cal O}\left(#1\right)}

\newcommand{\veeh}{V_h}
\newcommand{\ntrp}{{\mathcal I}_h}
\newcommand{\dubya}{W_h}
\newcommand{\Reyuls}{\mathbb{R}}
\newcommand{\Rn}{\Reyuls^n}

\newcommand{\beginproof}{\medskip\par\noindent{\bf Proof.\ }}
\newcommand{\proofend}{{\bf QED}\par\medskip}

\newcommand{\half}{{\textstyle{1\over 2}}}

\newcommand{\set}[2]{\left\lbrace #1 \; : \; #2 \right\rbrace}

\newcommand{\norm}[1]{\Vert #1 \Vert}

\newcommand{\tbnorm}[1]{\vert\kern-0.1em\vert\kern-0.1em\vert\, #1 \,\vert\kern-0.1em\vert\kern-0.1em\vert}

\newcommand{\qbnorm}[1]{\vert\kern-0.15em\vert\kern-0.15em\vert\kern-0.15em\vert\, #1 \,\vert\kern-0.15em\vert\kern-0.15em\vert\kern-0.15em\vert}

\newcommand{\teo}[1]{{\color{red} #1}}
\newcommand{\teoo}[1]{{\color{red} #1}}

\title{A Savitsky--Golay filter \\ for smoothing  Adaptive finite element computations}
\date{\today}

\author{
Teodoro Collin\thanks{
The University of Chicago, Chicago, Illinois \ \ 60636, USA; email:teoc@mit.edu};
current address, MIT \and
Gordon Kindlmann\thanks{
Department of Computer Science, The University of Chicago,
Chicago, Illinois \ \  60637, USA; email: glk@uchicago.edu} \and
L.~Ridgway Scott\thanks{
The University of Chicago,
Chicago, Illinois \ \  60637, USA; email: ridg@uchicago.edu}
}

\begin{document}

\maketitle



\begin{abstract}
The smoothing technique of Savitzky and Golay is extended to data defined on
multidimensional meshes.
\color{red}
We define a smoothness-increasing filter that is inspired by the
smoothness-increasing, accuracy-conserving (SIAC) filters.
Given the genealogy of the proposed smoother, we refer to it as
SISG, for smoothness-increasing, Savitzky--Golay.
SISG is \teoo{easily applied to finite-element computations on adaptive meshes.}
It can be applied effectively both for visualization and data analysis where the 
underlying approximations are discontinuous.
\color{black}
\end{abstract}

\section{Introduction}
\teo{
In finite element computations, many problems are best solved with non-smooth functions, but many post-processing tasks require continuous or smooth data.
Since post-processing tasks also require some sort of fidelity with
the original simulation data, post-processing tasks require smothers
that conserve the accuracy of simulations.
The smoothness increasing and accuracy conserving filters (SIAC)~\cite{lineSIACfinitevol,li2016smoothness,mirzargar2016smoothness,
  king2012smoothness}
have been proven to achieve these aims for a limited variety of geometries via a complicated convolution.
%
In a signal processing context, the Savitzky--Golay filter can achieve similar aims
on a regularly spaced sequence via a convolution that can be identified
with a collection of local least squares problems~\cite{savitzky1964smoothing}.
In this work, we merge the ideas of Savitzky and Golay and the SIAC methods to
develop a new smoother for finite element computations that requires limited
mesh assumptions and implementation effort.
Merging SIAC and Savitzky--Golay ideas, we call the method SISG.

In particular, SISG conserves accuracy and smooths data just like the SIAC filters,
but SISG smoothes via an orthogonal projection similar to that of Savitzky--Golay.
Since SISG is based on $L^2$ projection and since many automated finite
elements systems (such as FEniCS~\cite{FEniCSbook} or Firedrake~\cite{rathgeber2017firedrake})
readily provide many useful $L^2$ projections, SISG is easy to implement.
Further, SISG relies on typical mesh geometry assumptions so SISG is especially
applicable to highly refined and irregular meshes
such as those arising from adaptive methods~\cite{cockburn2003enhanced,king2012smoothness}.
In contrast, SIAC requires non-trivial implementation effort
beyond a typical finite element system and, although there
are promising empirical studies of adaptive versions of
SIAC~\cite{jallepalli2019adaptive}, SIAC has only been
verified for irregular meshes in certain
circumstances~\cite{li2016smoothness}.
The purpose of this paper is to provide the theoretical basis
for SISG and to explore it more critically and genealogically.}

The approach of Savitzky and Golay~\cite{savitzky1964smoothing,lrsBIBbq}
to data smoothing involves fitting data via least squares to a polynomial
in a window of fixed size.
It is applied to data values $y_i$ corresponding to time (or other variable)
points $x_i$ that are spaced uniformly, that is, $x_i=i h$ for some $h>0$.
Using the method of least squares, a polynomial $P$ of degree $k$ is chosen
so that the expression
$$
\sum_{i=1}^{r} (y_i-P(x_i))^2
$$
is minimized over all polynomials of degree $k$.
A major result of \cite{savitzky1964smoothing} was the identification of the
least squares processes and subsequent evaluation of approximations as being
equivalent to a convolution with a discrete kernel of finite extent
in each case.
\teo{An even more efficient algorithm in the cubic case was explored in \cite{lrsBIBbq}.}
Least-squares approaches in the context of signal processing
was compared with interpolation in~\cite{unser1997approximation}
with a special emphasis on theoretical tools developed
originally for the finite element method \cite{strang2011}.
See also \cite{van2002least}.
The polynomial $P$ can be used in various ways: providing an approximation to the
data at some point in the window $x_{1},\dots,x_{r}$, or a derivative
at such a point, or a second derivative, and so forth.
%

The SIAC filters derive from the work of Bramble and Schatz \cite{bramble1977higher}
and Thom{\' e}e \cite{thomee1977high} who realized that higher-order local accuracy
could be extracted by averaging finite element approximations.
\teo{In particular, Thom{\' e}e \cite{thomee1977high} proved his results via
examining the oscillation in the error and the impact of B-spline kernels
on this error oscillation.}
This was further developed in the context of finite element approximations
of hyperbolic equations \cite{cockburn2003enhanced}.
It was also realized that this approach could be extended to compute accurate
approximations to derivatives even when the original finite elements are
discontinuous \cite{ryan2009local}.
Thus, the SIAC (smoothness-increasing accuracy-conserving) filters
\cite{lineSIACfinitevol,li2016smoothness} were developed
\begin{enumerate}
\item
to smooth discontinuities in functions or their derivatives and
\item
to extract extra accuracy associated with error oscillation in certain settings.
\end{enumerate}

\teo{SIAC achieves its goals via convolution with a localized kernel, generated
via B-Splines.
When meshes are irregular, this can pose certain difficulties
as the relationship between the extents of mesh elements
and the design of the kernels is critical~\cite{li2016smoothness,jallepalli2019adaptive,Li2019}.
Moreover, many problems lacking solution regularity simultaneously feature limited error oscillations
and highly irregular meshes used to resolve the solution
to high accuracy~\cite{lrsBIBih}.
}

\begin{figure}
\centerline{(a)\includegraphics[width=3.0in]{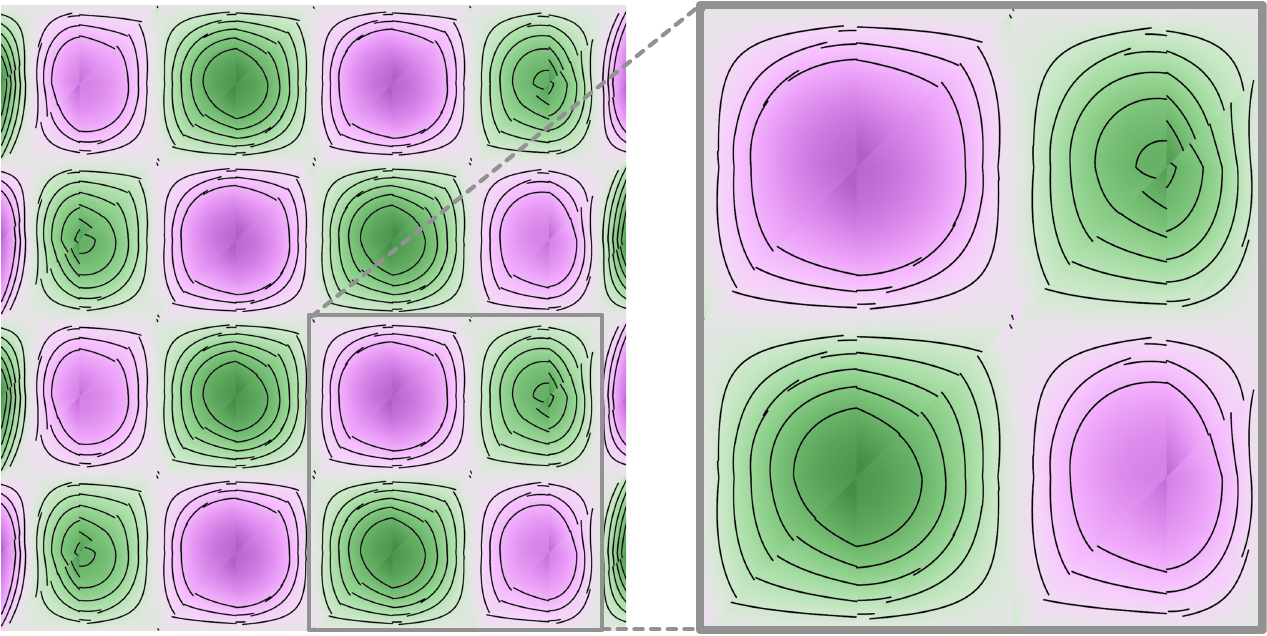}
\quad      (b)\includegraphics[width=3.0in]{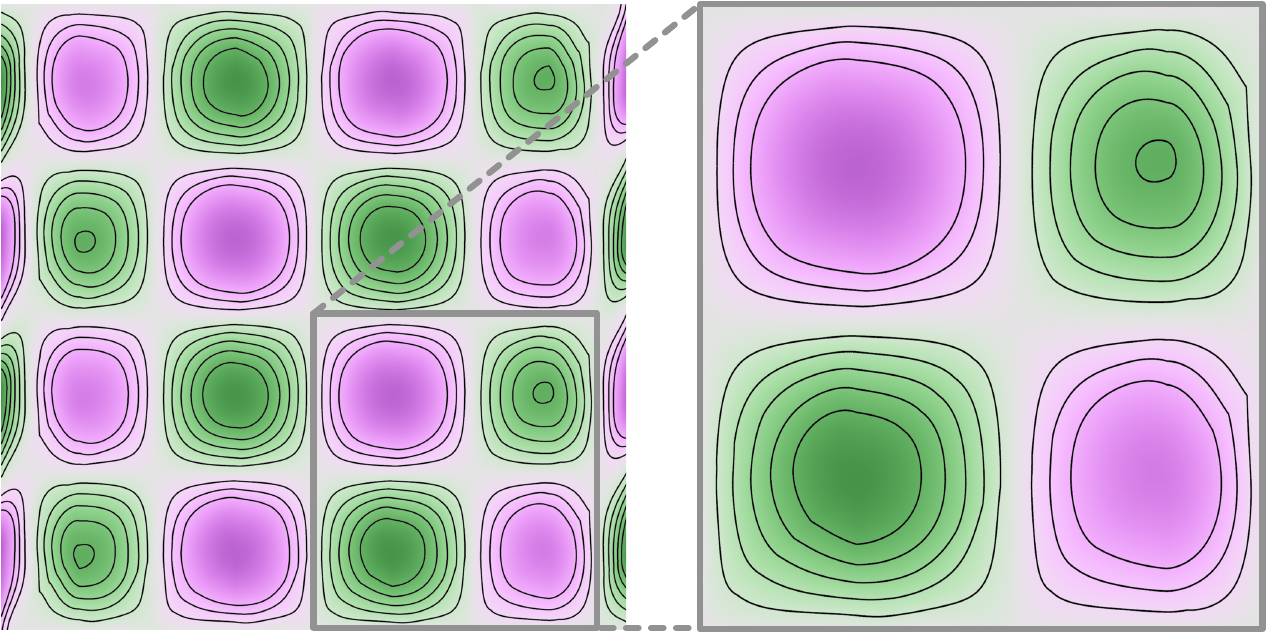}}
\caption{$x$-derivative of solution: (a) before filtering, (b) after filtering. A zoomed view is provided to facilitate comparisons within and between figures.}
\label{fig:bothsiac}
\end{figure}

In this work, we utilize the original idea of Savitzky and Golay \cite{savitzky1964smoothing}
to smooth or differentiate data defined via a finite element representation
on a multidimensional mesh $M_h$, where $h$ measures in some way the size
of the mesh elements.
Our method achieves goal 1.~of SIAC filters, without any assumption of error oscillation.
This is illustrated in Figure \ref{fig:bothsiac} where we apply our approach
to smooth the derivative of a standard finite-element approximation.
Our filter does not appear to achieve goal 2., but it does apply with highly
refined meshes.
It is ``accuracy conserving'' in a certain sense, so the term SIAC still is
applicable, but it is not ``accuracy improving'' in the sense of
\cite{bramble1977higher,thomee1977high,cockburn2003enhanced}.
Thus we introduce a different acronym for the technique
introduced here, SISG, for smoothness increasing Savitzky--Golay.

\teo{Although we have borrowed from the SIAC and Savitzky--Golay ideas,
what we propose is fundamentally different and has a different application
domain.
Our domain is highly refined meshes.
The original SIAC and Savitzky--Golay techniques are closely associated with regular meshes.
Further, both techniques are also computed with convolution.
On a highly refined mesh, different scales are required, as opposed to convolution
with a fixed mesh size.
Moreover, convolution on a nonuniform mesh requires
the careful and costly construction of a regular mesh
overlayed on the nonuniform mesh~\cite{jallepalli2019adaptive}.
The method we propose utilizes a global projection scheme for efficiency, ease of use, and elegance.
Another distinction from the Savitzky--Golay
method is that our approach works at the level
of functions rather than data points.
Given a nonsmooth function, we produce smoother functions.
This has several advantages with regard to visualization
and data extraction techniques.

Additionally, although we focus on merging SIAC and Savitzky--Golay in this work,
we suspect other connections may emerge.
For instance, the Savitzky--Golay approach of using discontinuous polynomials
for local approximation to define smoothed data is similar to what is done
in Cl{\'e}ment Interpolation \cite{ref:ClementInterpolation}.
Similarly, SIAC has recently been connected to the theory of
quasi-interpolation, which then relates in a number of ways
to splines used to generate other convolution-based, least-squares processes \cite{deBoor1993,unser1997approximation}.
We also provide a similar disclaimer for our choice of applications;
we will focus on visualization as an application, but smoothing has
many applications.
For instance, an intriguing application lies in residual-based error
estimators~\cite{dedner2019residual}.
}

\subsection{Sobolev spaces}

To quantify the notions of smoothness and accuracy in SIAC, Sobolev spaces
are used, defined as follows.
Let $\Omega\subset\Rn$ be a domain of interest, and define
$$
\norm{u}_{L^2(\Omega)}=\bigg(\int_\Omega |u(x)|^2\,dx\bigg)^{1/2}
$$
for a scalar or vector-valued function $u$.
Let $\nabla^m$ denote the tensor of $m$-th order partial derivatives.
Then
$$
\norm{u}_{s}=\sum_{m=0}^s \norm{\,|\nabla^m u|\,}_{L^2(\Omega)},
$$
where $|T|$ denotes the Euclidean norm of a tensor $T$ of arity $\alpha$,
thought of as a vector of dimension $n^\alpha$.
Define $H^s(\Omega)$ to be the set
of functions whose derivatives of order up to $s$ are square integrable,
that is, such that $\norm{u}_{s}<\infty$.
Note that $\norm{u}_{0}=\norm{u}_{L^2(\Omega)}$.
More generally, we write
$$
\norm{T}_{s}=\sum_{m=0}^s \norm{\,|\nabla^m T|\,}_{L^2(\Omega)},
$$
for a tensor-valued function $T$.
Note that for a tensor $T$ of arity $\alpha$, $\nabla^m T$ is
a tensor of arity $\alpha+m$.
These norms can be extended to allow non-integer values of $s$ \cite{lrsBIBgd}.

\subsection{Approximation on meshes}

Suppose that $\Omega\subset\Rn$ is subdivided in some way by a mesh $M_h$ (e.g., a 
triangulation, quadrilaterals, prisms, etc.). 
Suppose further that $W_h$ is a finite element space defined on $M_h$ and
that $u_h\in\dubya$ is some approximation to a function $u\in H^s(\Omega)$.
In many cases \cite{lrsBIBgd}, an error estimate holds of the form
\color{red}
\begin{equation}\label{eqn:neuassum}
\norm{u- u_h}_{L^2(\Omega)}^2 \leq C \sum_e h_e^{2(k+1)}\norm{\nabla^{k+1}u}_{L^2(e)}^2.
\end{equation}
\color{black}
The approximation $u_h$ could be defined by many techniques, including Galerkin
approximations to solutions of partial differential equations \cite{li2016smoothness}.
\color{red}
The approximation $u_h$ can be a vector, or tensor, and it could for example come
from taking the gradient of a finite-element approximation.
\color{black}
The SIAC objective is to create an operator $\Pi_h$ that maps $u_h$ into
a smoother space $\veeh$ in a way that maintains this accuracy:
\begin{equation}\label{eqn:getumd}
\norm{u-\Pi_h u_h}_{-t}\leq C h^{s+t}\norm{u-u_h}_{s}, \; 0\leq t\leq t_0,
\end{equation}
\color{red}
where $\norm{v}_{s}$ denotes the norm in $H^s(\Omega)$,
 and when $s$ is negative, the norm in $H^s$ is defined by duality \cite{lrsBIBgd}.
\color{black}
Moreover, provided that $\veeh\subset H^\tau(\Omega)$ for $\tau>0$, we will
also be able to show that 
\begin{equation}\label{eqn:derivd}
\norm{u-\Pi_h u_h}_{\tau}\leq C h^{s-\tau}\norm{u-u_h}_{s}.
\end{equation}
This means that the SISG filter $\Pi_h u_h$ can provide optimal-order approximations 
of derivatives of $u$, even though the space $W_h$ from which $u_h$ comes 
may harbor discontinuous functions.

We will define $\Pi_h$ as the $L^2(\Omega)$ projection of $W_h$ onto $\veeh$.
In this way, our proposed SIAC is very similar to the unified Stokes algorithm (USA) proposed in \cite{lrsBIBia}.
%

\section{Savitzky-Golay as a projection}
\label{sec:sgm}
\teo{In order to see a Savitzky-Golay filter as a predecessor to our SISG, we must carefully examine the filter. Savitzky-Golay presumes a  finite sequence, $\{x_{i}\}_{i=0}^{M}$, with a spacing parameter $h > 0$ that represents the spacing of sample locations. Savitzky-Golay aims to produce a smoother sequence $\{y_{i}\}_{i=0}^{M}$. Savitzky-Golay produces the sequence $\{y_{i}\}$ via evaluating polynomials produced via a sequence of local projections on particular subsequences of $\{x_{i}\}$ of length $r$ in an inner product that we will now define.

First, given two sequences of real numbers $f$ and $g$ of length $r$, we define the inner product

\begin{equation}\label{eqn:sgipmd}
(f,g)_{r}= \sum_{i=0}^{r} f_i g_i.
\end{equation}

Then we extend this inner product to polynomials of degree $k$ via

\begin{equation}\label{eqn:polymd}
(P,Q)_{r,j}= \sum_{i=0}^{r} P((i+j)h) Q((i+j)h)
\end{equation}
with $j > 0$.
More precisely, the inner-products are defined on $\Reyuls^{r}$, and we think of the space of polynomials of degree $k$ as a $k+1$ dimensional subspace via the identification $f_i=P((i+j)h)$ for all $0\leq i\leq r$.
To be an inner-product on polynomials of degree $k$, we must have $r>k$, for
otherwise a nonzero polynomial of degree $k$ could vanish at all of the grid points.

The desired inner product then operates so that we can solve the following least squares problem: for a fixed sequence $\{z_{i}\}$ of length $r$ and an integer $j$, find the polynomial $P$ that minimizes
\begin{equation}\label{eqn:asinpd}
(P-z,P-z)_{r,j} = \sum_{i=0}^{r} (P((i+j)h) - z_i)^2.
\end{equation}

Savitzky-Golay uses $(\cdot, \cdot)_{r,j}$ to produce $\{y_{i}\}$ from $\{x_{i}\}$ as follows. For each sub-sequence of $\{x_{i}\}$ of the form $\{x_{j+s}\}_{s=0}^{s=r} = \{z_{j,s}\}$, Savitzky-Golay uses least squares methods to find $P_{j}$ that minimizes $(P_{j}-z_{j}, P_{j}-z_{j})_{r,j}$. With the resulting sequence, $\{P_{j}\}$, Savitzky-Golay produces $y_{j}= P_{j}(c_{j})$ where $c_{j}$ is the median value of $\{(s+j)h\}_{s=0}^{r}$.

Savitzky-Golay is efficient because the least squares process given by~\eqref{eqn:asinpd} can be identified with a shift invariant projection, meaning the simultaneous application of the projection across all $j$ can be identified as a convolution~\cite{Persson2003}.
Indeed, Savitzky-Golay sometimes precomputes a projection via a QR factorization and then uses this projection to produce convolution coefficients~\cite{williampress1992}.

However, although Savitzky-Golay is efficient because it can be a convolution, Savitzky-Golay was designed as a projection~\cite{savitzky1964smoothing}.
Indeed, we can see Savitzky-Golay as a single projection onto a single global space via turning each Savitzky-Golay window into a local inner product space $(\Reyuls^r, (\cdot, \cdot)_{r,j})$ for a window or cell $j$ and then taking the cartesian product of these spaces over $j$.
We can define a least squares problem on this global space and this problem would be equivalent to solving~\eqref{eqn:asinpd} for all $j$.
Though it would be comparatively inefficient and a bit of a hack, we believe that many finite element packages could implement Savitzky-Golay in this manner.

Therefore, a Savitzky-Golay filtering could be applied to other inputs via a similar aggregation of related local projections. Savitzky-Golay has been viewed in this manner in order to create other filters in signal processing~\cite{GSG} and we will view it in this manner to produce smoothness increasing filters that unlike the previously mentioned SIAC filters use readily accessible elements of both finite element theory and implementation.
}


\section{Least-squares projections for finite elements}
\label{sec:lsp}
We can define the projection $\Pi_h$ from $L^2(\Omega)$ to $\veeh$ by
\begin{equation}\label{eqn:defpih}
(\Pi_h u,v)_{L^2(\Omega)}= (u,v)_{L^2(\Omega)}\quad\forall v\in\veeh.
\end{equation}
This finite dimensional system of equations is easily defined and solved
in automated systems \cite{rathgeber2017firedrake,lrsBIBih}.
Although this requires the solution of a global system, it allows the use of
general meshes, including highly graded meshes.
It is convenient to restrict to ones that are nondegenerate, by which 
we mean that each element $e$ of the mesh (simplex, cube, prism, etc.) 
has the following properties.

\subsection{Mesh assumptions}
\label{sec:ma}
\begin{definition}\label{def:newnondeg}
Let $D$ be a finite, positive integer.
Let $M_h$ be a family of meshes consisting of elements $e\in M_h$.
For each $e\in M_h$, let $\rho_e$ be the diameter of the largest sphere 
contained in $e$, and let $h_e$ be the diameter of the smallest sphere 
containing $e$.
We say that the family of meshes $M_h$ is {\bf nondegenerate} if
there is a constant $\gamma$ such that, for all $h$ and $e$,
\begin{equation}\label{eqn:newnondeg}
h_e\leq \gamma \rho_e
\end{equation}
and such that $e$ is a union of at most $D$ domains that are
star-shaped \cite{lrsBIBat} with respect to a ball of radius $h_e/\gamma$.
\end{definition}

In most cases \cite{lrsBIBgd}, $D=1$, that is, the elements are star-shaped
with respect to a ball.
This holds, for example, if the elements are all convex.
Thus the usual definition \cite{lrsBIBgd} is based just on \eqref{eqn:newnondeg}.
But we have allowed $D>1$ for generality.

Note that the parameter $h$ is not a single mesh size but rather a function
defined on the mesh that prescribes the local mesh size.
If the quantifiers guarding \eqref{eqn:newnondeg} are rearranged, so that
\begin{equation}\label{eqn:newquasiu}
\max_e h_e\leq \gamma \min_e \rho_e,
\end{equation}
then we call the mesh quasi-uniform \cite{lrsBIBgd}, and we can define a single
mesh size $h=\max_e h_e$.

\subsection{Approximation theory}

For a nondegenerate family of meshes $M_h$, there is a constant $C$ depending only
on $\gamma$ and $D$, such that, 
for all $e$ and $h$, there is a polynomial $P_e$ of degree $k$ with the property
\begin{equation}\label{eqn:dgonlihyp}
\norm{\nabla^m(u-P_e)}_{L^2(e)} \leq C  h_e^{k+1-m}
\norm{\nabla^{k+1} u}_{L^2(e)}.
\end{equation}
This is a consequence of the Bramble-Hilbert lemma \cite{lrsBIBat}.
It in particular demonstrates adaptive approximation using discontinuous
Galerkin methods.

An estimate similar to \eqref{eqn:dgonlihyp}
involving an approximation operator $\ntrp$, called an interpolant,
is well known for essentially all finite-element methods.
Thus we make the following {\em hypothesis}:
\begin{equation}\label{eqn:hyppih}
\sum_{m=0}^t \sum_e h_e^{2m}\norm{\nabla^m(u-\ntrp u)}_{L^2(e)}^2
\leq C_I \sum_e h_e^{2(k+1)} \norm{\nabla^{k+1} u}_{L^2(e)}^2.
\end{equation}
\color{red}
The constant $C_I$ typically depends only on $\gamma$ and $D$ in 
Definition \ref{def:newnondeg}, but for generality we do not assume this.
\color{black}
Note that we have apportioned the quantity $h_e^{k+1-m}$ in \eqref{eqn:dgonlihyp}
partly on the left and partly on the right in \eqref{eqn:hyppih}.
This balancing act is arbitrary and could be done in a different way.
But the thinking in \eqref{eqn:hyppih} is that we have chosen the mesh so
that the terms on the right-hand side of the inequality in \eqref{eqn:hyppih}
are balanced in some way.
That is, we make the mesh size small where the derivatives of $u$ are large.

For discontinuous Galerkin methods, we take $\ntrp u|_e =P_e$ on each 
element $e\in M_h$, where $P_e$ comes from \eqref{eqn:dgonlihyp}.
On a quasi-uniform \cite{lrsBIBgd} mesh of size $h$, \eqref{eqn:hyppih} simplifies to
$$
\sum_{m=0}^t h^{m}\norm{\nabla^m(u-\ntrp u)}_0 \leq C h^{k+1} \norm{u}_{k+1}.
$$

The following can be found in \cite{lrsBIBgd} and other sources.

\begin{lemma}\label{lem:hyppih}
Let $\veeh$ consist of piecewise polynomials that include complete polynomials
of degree $k$ on each element of a nondegenerate subdivision of $\Omega$, such as 
\begin{itemize}
\item
($t=0$) Discontinuous Galerkin (DG) ($k\geq 0$), 
\item
($t=1$) Lagrange ($k\geq 1$), Hermite ($k\geq 3$), tensor-product elements ($k\geq 1$),
\item
($t=2$) Argyris ($k\geq 5$ in two dimensions),
\end{itemize}
as well as many others.
Then the hypothesis \eqref{eqn:hyppih} holds with $t$ and $k$ as specified.
\end{lemma}

\subsection{Estimates for the projection}

\color{red}
We first show that the SISG filter provides good approximation to smooth functions.
\color{black}

\begin{theorem}\label{thm:hyppih}
Under the hypothesis \eqref{eqn:hyppih}, and in particular for the spaces
$\veeh$ listed in Lemma \ref{lem:hyppih}, there is a constant $C$ such that
\begin{equation}\label{eqn:extrahyppih}
\sum_{m=0}^t \sum_e h_e^{2m}\norm{\nabla^m(u-\Pi_h u)}_{L^2(e)}^2
\leq C  \sum_e h_e^{2(k+1)} \norm{\nabla^{k+1}u}_{L^2(e)}^2
\end{equation}
for any nondegenerate family of meshes.
The constant $C$ depends only on $\gamma$ and $D$ in Definition \ref{def:newnondeg}
\color{red}
and the constant $C_I$ in \eqref{eqn:hyppih}.
\color{black}
\end{theorem}

\beginproof
First of all, consider the case $t=0$.
In this case, the result to be proved is
\begin{equation}\label{eqn:extoporvpih}
\sum_e \norm{u-\Pi_h u}_{L^2(e)}^2
\leq C  \sum_e h_e^{2(k+1)} \norm{\nabla^{k+1}u}_{L^2(e)}^2.
\end{equation}
Since the $L^2$ projection provides optimal approximation in $L^2$, we find
$$
\sum_e \norm{u-\Pi_h u}_{L^2(e)}^2=\norm{u-\Pi_h u}_{L^2(\Omega)}^2
\leq \norm{u-\ntrp u}_{L^2(\Omega)}^2.
$$
Then the hypothesis \eqref{eqn:hyppih} implies \eqref{eqn:extoporvpih}, which
completes the proof in the case $t=0$.

Using the triangle inequality and standard inverse estimates \cite{lrsBIBgd}
for nondegenerate meshes, we have
\begin{equation}\label{eqn:nvrstopih}
\begin{split}
\norm{\nabla^m(u-\Pi_h u)}_{L^2(e)}&\leq 
\norm{\nabla^m(u-\ntrp u)}_{L^2(e)}+ \norm{\nabla^m(\ntrp u-\Pi_h u)}_{L^2(e)} \\
&\leq \norm{\nabla^m(u-\ntrp u)}_{L^2(e)}+ C h_e^{-m} \norm{\ntrp u-\Pi_h u}_{L^2(e)} \\
&\leq \norm{\nabla^m(u-\ntrp u)}_{L^2(e)}+ C h_e^{-m}\big( \norm{u-\ntrp u}_{L^2(e)} 
+  \norm{u-\Pi_h u}_{L^2(e)}\big) .
\end{split}
\end{equation}
Squaring and summing this over $e$, and using the hypothesis \eqref{eqn:hyppih} yields
$$
 \sum_e h_e^{2m}\norm{\nabla^m(u-\Pi_h u)}_{L^2(e)}^2
\leq C \Big(\sum_e h_e^{2(k+1)} \norm{\nabla^{k+1}u}_{L^2(e)}^2
+\norm{u-\Pi_h u}_{L^2(e)}^2\Big),
$$
with a possibly larger constant $C$.
Using \eqref{eqn:extoporvpih} and summing over $m$ completes the proof.
\proofend

\color{red}
\subsection{Estimates for the SISG filter}
\color{black}

The following theorem justifies the 
\color{red}
SI part of SISG for $\Pi_h$.
\color{black}

\begin{theorem}\label{thm:neuassum}
\color{red}
Suppose that hypothesis \eqref{eqn:hyppih} holds.
\color{black}
Suppose $\veeh$ is a space as in Lemma \ref{lem:hyppih}.
Under 
the assumption \eqref{eqn:neuassum}, we have
\begin{equation}\label{eqn:neuthmrslt}
\sum_{m=0}^t \sum_e h_e^{2m}\norm{\nabla^m(u-\Pi_h u_h)}_{L^2(e)}^2
\leq C  \sum_e h_e^{2(k+1)} \norm{\nabla^{k+1}u}_{L^2(e)}^2.
\end{equation}
\color{red}
The constant $C$ depends only on $\gamma$ and $D$ in Definition \ref{def:newnondeg}
and the constant $C_I$ in \eqref{eqn:hyppih}.
\color{black}
\end{theorem}

\beginproof
The proof is similar to that of Theorem \ref{thm:hyppih}.
We write 
$$
u-\Pi_h u_h= u-\ntrp u+\ntrp u-\Pi_h u_h. 
$$
The required estimate for $u-\ntrp u$ is hypothesis \eqref{eqn:hyppih}.
Using inverse estimates \cite{lrsBIBgd} for nondegenerate meshes, we find
\begin{equation}\label{eqn:neurstopih}
\begin{split}
 h_e^{-m} \norm{\nabla^m(\ntrp u-\Pi_h u_h)}_{L^2(e)}
&\leq  C \norm{\ntrp u-\Pi_h u_h}_{L^2(e)} \\
&\leq  C \big( \norm{u-\ntrp u}_{L^2(e)} +  \norm{u-\Pi_h u_h}_{L^2(e)}\big) .
\end{split}
\end{equation}
Squaring and summing over $e$, and applying \eqref{eqn:hyppih},
the triangle inequality, and \eqref{eqn:extrahyppih}, yields
\begin{equation}\label{eqn:sqsumtopih}
\begin{split}
\sum_e h_e^{-2m} &\norm{\nabla^m(\ntrp u-\Pi_h u_h)}_{L^2(e)}^2
\leq  2C \Big(\sum_e \norm{u-\ntrp u}_{L^2(e)}^2 
             + \sum_e \norm{u-\Pi_h u_h}_{L^2(e)}^2\Big) \\
&\leq C' \Big( \sum_e h_e^{2(k+1)} \norm{\nabla^{k+1}u}_{L^2(e)}^2
             + \sum_e \norm{u-\Pi_h u}_{L^2(e)}^2
             + \sum_e \norm{\Pi_h (u-u_h)}_{L^2(e)}^2\Big), \\
&\leq C'' \Big( \sum_e h_e^{2(k+1)} \norm{\nabla^{k+1}u}_{L^2(e)}^2
             + \sum_e \norm{\Pi_h (u-u_h)}_{L^2(e)}^2\Big), \\
\end{split}
\end{equation}
for appropriate constants $C'$ and $C''$.
Since $\Pi_h$ is the $L^2(\Omega)$ projection,
$$
 \sum_e \norm{\Pi_h (u-u_h)}_{L^2(e)}^2=
\norm{\Pi_h (u-u_h)}_{L^2(\Omega)}^2\leq \norm{u-u_h}_{L^2(\Omega)}^2.
$$
Thus applying \eqref{eqn:neuassum} completes the proof.
\proofend

\begin{figure}
\centerline{(a)\includegraphics[width=3.0in]{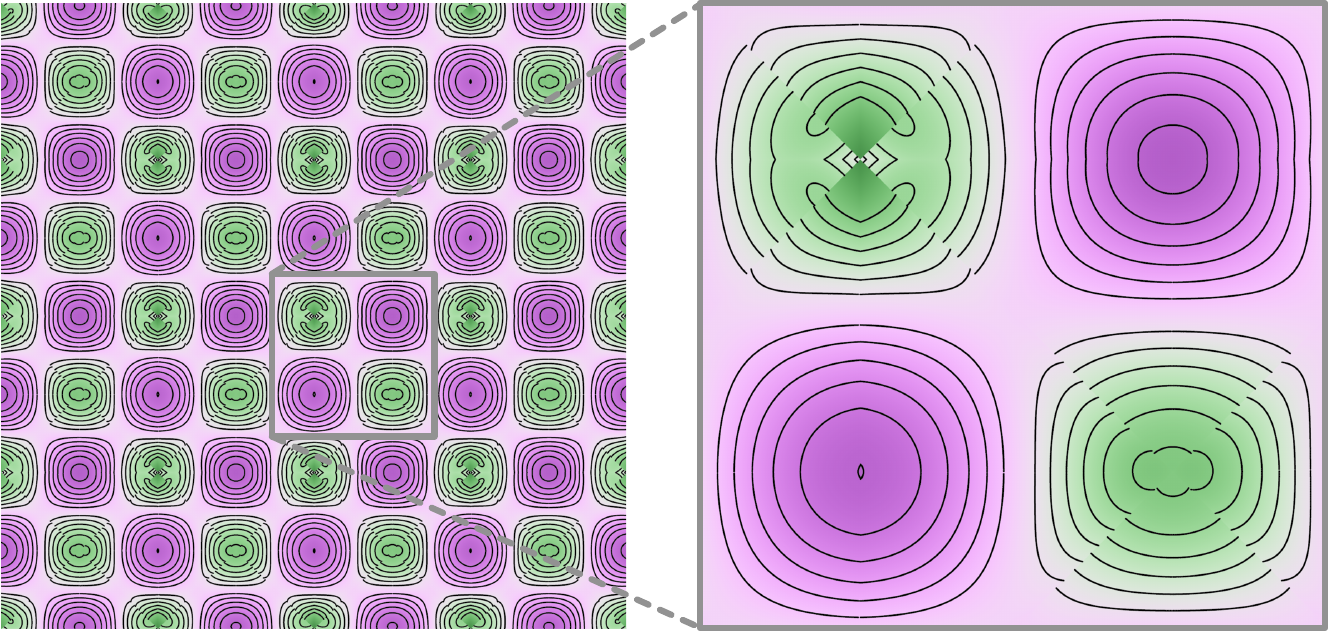}
\quad      (b)\includegraphics[width=3.0in]{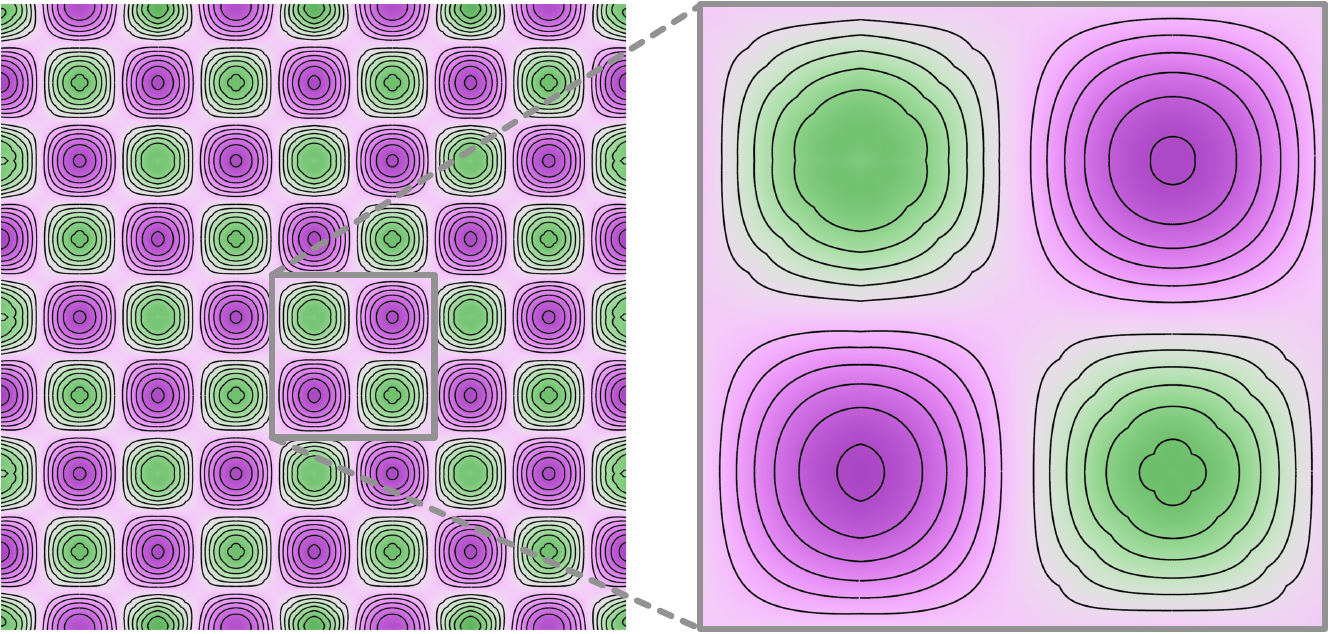}}
\caption{Pressure in the solution of the Stokes equation:
(a) standard Scott-Vogelius discontinuous pressure, (b) after filtering using
SISG. A zoomed view is provided to facilitate comparisons within and between figures.}
\label{fig:stokex}
\end{figure}
\teo{
\section{Theoretical and Implementation comparison with SIAC}
In this section, we will briefly compare the theoretical and implementation aspects of SIAC and SISG. 
In section \ref{sec:hirefme}, we give performance numbers for SISG on a highly refined mesh. Such problems  lack the domain regularity necessary for the improved negative norm accuracy which underpins the AC aspect of SIAC filters.

As we have shown in the previous section, the key theoretical assumptions of SISG are staples of finite element theory~\cite{lrsBIBgd}. In particular, the mesh geometry assumptions of Sec.~\ref{sec:ma} apply to many practical meshes, including many found in adaptive settings, whereas the story for SIAC is more complicated. Originally, SIAC was developed for uniform, translation invariant, or structured triangular meshes, but SIAC has and is being generalized to other geometric situations; see \cite{jallepalli2019adaptive} for a recent overview. SIAC works under the following mesh geometry condition:

\begin{definition}
\label{def:kings}
Let $M\subset\mathbb{R}^{n}$ be a mesh. For each element, $e\in M$, define $h_{e}$ to be the $n-$tuple of real numbers such that the $j$-th element is the maximum extent of $e$ in the $j$-th coordinate direction. Then we say that $M$ has an integer partitioning property if there is a fixed real number $H > 0$ and there exists for every $e\in M$ an integer tuple $l_{e}$ such that $h_{e,i} = H/l_{e,i}$ for $1\leq i\leq n$.
\end{definition}

Definition \ref{def:kings} corresponds to Theorem 3.1 in \cite{king2012smoothness} and can be more easily understood by examining Figure 5 in that paper. Definition \ref{def:kings} is the most general proven theoretical condition that we are aware of and it falls quite far from the generality found in many meshes covered by standard assumptions such as those in Definition \ref{def:newnondeg}. Outside of this condition, SIAC techniques have also been extended to meshes with periodic boundary conditions or to meshes in 1D~\cite{li2016smoothness,lithesis}.

Beyond theory, SIAC has also been studied empirically and several studies suggest more general results via various means. For example, \cite{jallepalli2019adaptive} proposes locally varying the characteristic length based on a simple indicator whereas \cite{Curtis2008} proposes using the local $L^{2}$ projection to transfer solutions on smoothly varying meshes (affine transformations of uniform meshes) to a locally uniform mesh. We can't comment on the theoretical results that these empirical studies might imply, but we can analyze the implementation aspects of SIAC and SIAC variants.

Unlike SISG, SIAC and SIAC variants (e.g. L-SIAC) require another layer of implementation that most FEM systems do not typically supply. In particular, both SIAC and L-SIAC require a convolution integral to be subdivided to avoid discontinuities created by both mesh element boundaries and the involved B-Spline kernels. In the L-SIAC case, the convolution integral is over a parametric curve that must subdivided at intersections with mesh element boundaries~\cite{jallepalli2019adaptive, Jallepalli2019}. Variants of SIAC disregard breaks caused by the B-Splines kernels, but these incur numerical crimes, meaning they are tools for specific computational situations (e.g. when subdividing as much as necessary is impossible)~\cite{Mirzaee2010}. To our knowledge, these techniques don't fall out of existing FEM software in the same way that SISG or other techniques do. For instance, a potential L-SIAC pre-processing step for computing an adaptive characteristics length could be computed via tricks in UFL (e.g. using cell diameters and vertices~\cite{Alns2014}), but the (possibly implicit) creation of a new mesh for convolution integrals will require software that does not currently ship with many FEM packages. 

}
\section{Applications}

There are situations in which a finite-element approximation, or its derivative,
is naturally discontinuous.
The SIAC-like operator proposed here is then essential if we want 
to visualize the corresponding approximations accurately.
We consider three examples in which the choice of the approximation space
$W_h$ is dictated by the structure of the problem.
Switching to a smoother space may give suboptimal results, so we are forced
to deal with visualization of discontinuous objects.

\subsection{Derivatives of standard Galerkin approximations}

Many standard finite element approximations \cite{lrsBIBgd,lrsBIBih}
use continuous (but not $C^1$) finite elements.
It is only recently \cite{kirby2018code} that $C^1$ elements have become 
available in automated finite-element systems.
Thus visualization of derivatives of finite-element approximations requires
dealing with discontinuous piecewise polynomials.

\color{red}
The Galerkin method optimizes the approximation of the gradient $\nabla v$
of the solution $v$ of a partial differential equation, not the function itself.
Thus we can think of the finite-element approximation $u_h=\nabla v_h$ as
the primary variable in the optimization.
The approximation $u_h$ of $u=\nabla v$ satifies the estimate \eqref{eqn:neuassum}
because of Ce{\' a}'s Lemma \cite{lrsBIBgd} which guarantees quasi-optimal
approximation of the gradients.
In many cases, better estimates of function value errors do not hold.
Thus we consider this as the base case for testing SISG.
\color{black}

We take as an example the equation
$$
-\Delta v = 32 \pi^2 \cos(4\pi x) \sin(4\pi y)\quad\hbox{in}\;[0,1]^2,
\qquad v=0\quad\hbox{on}\;\partial[0,1]^2,
$$
which has the exact solution $v(x,y)=-\cos(4\pi x) \sin(4\pi y)$.
This was approximated on an $8\times 8$ regular mesh of squares divided
into two right triangles, using continuous piecewise polynomials
of degree 3, the resulting approximation being denoted by $v_h$.
Depicted in Figure \ref{fig:bothsiac} is the (discontinuous) derivative 
$\partial_x v_h$, together with the smoothed version $\Pi_h(\partial_x v_h)$.
\color{red}
This demonstrates the power of SISG to smooth the singular $x$-component of $u_h$.
\color{black}

\subsection{Visualizing the pressure in Stokes}

Most finite-element methods for solving the Stokes, Navier-Stokes, or
non-Newtonian flow equations \cite{lrsBIBih,lrsBIBij} involve a 
discontinuous approximation of the pressure.
The unified Stokes algorithm (USA) proposed in \cite{lrsBIBia} uses 
the projection method described here to smooth the pressure.
In Figure \ref{fig:stokex}, we present the example from \cite{lrsBIBia}.
In this case, $\Pi_h p_h$ is exactly the USA pressure.

\subsection{DG for mixed methods}

\begin{figure}
\centerline{(a)\includegraphics[width=3.0in]{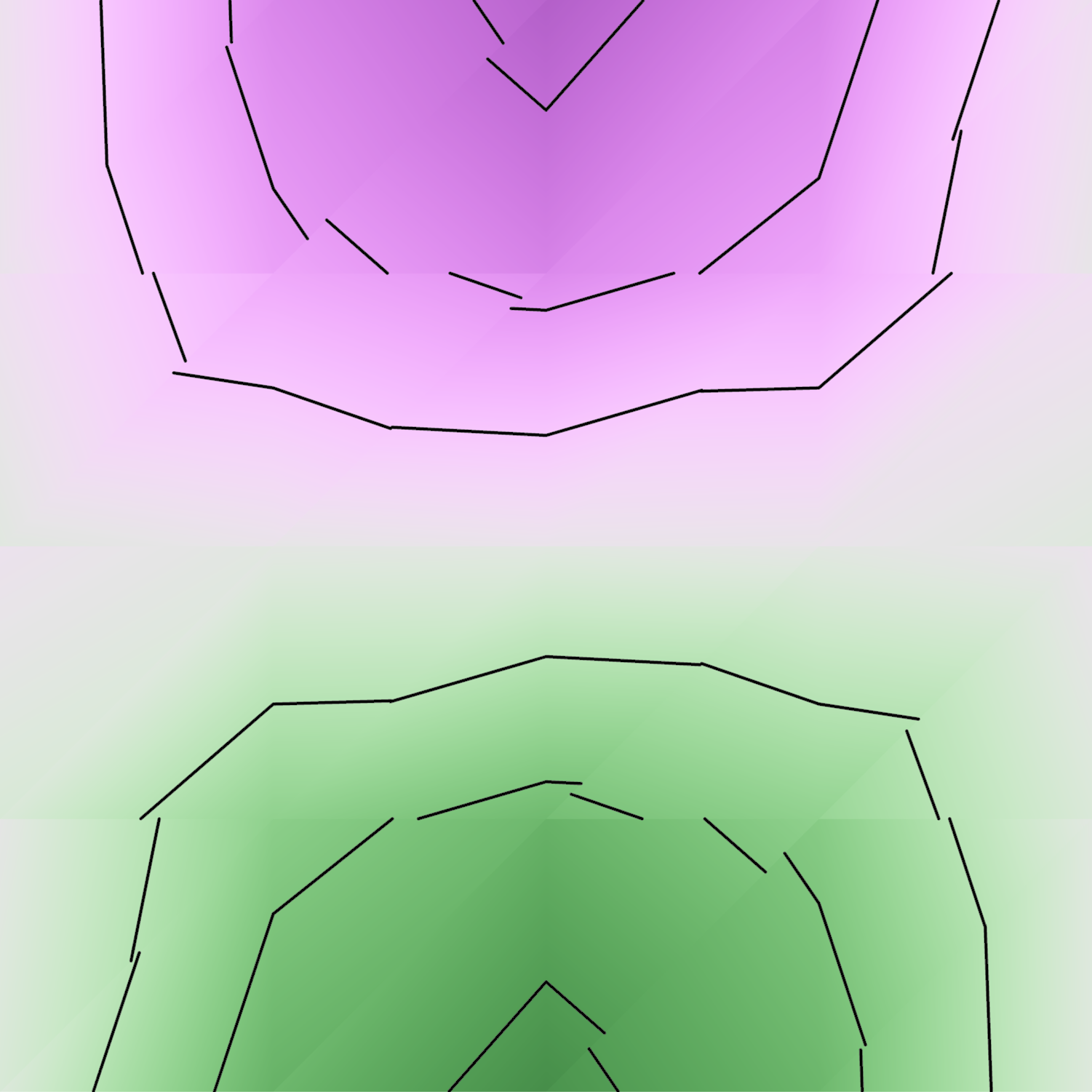}
\quad      (b)\includegraphics[width=3.0in]{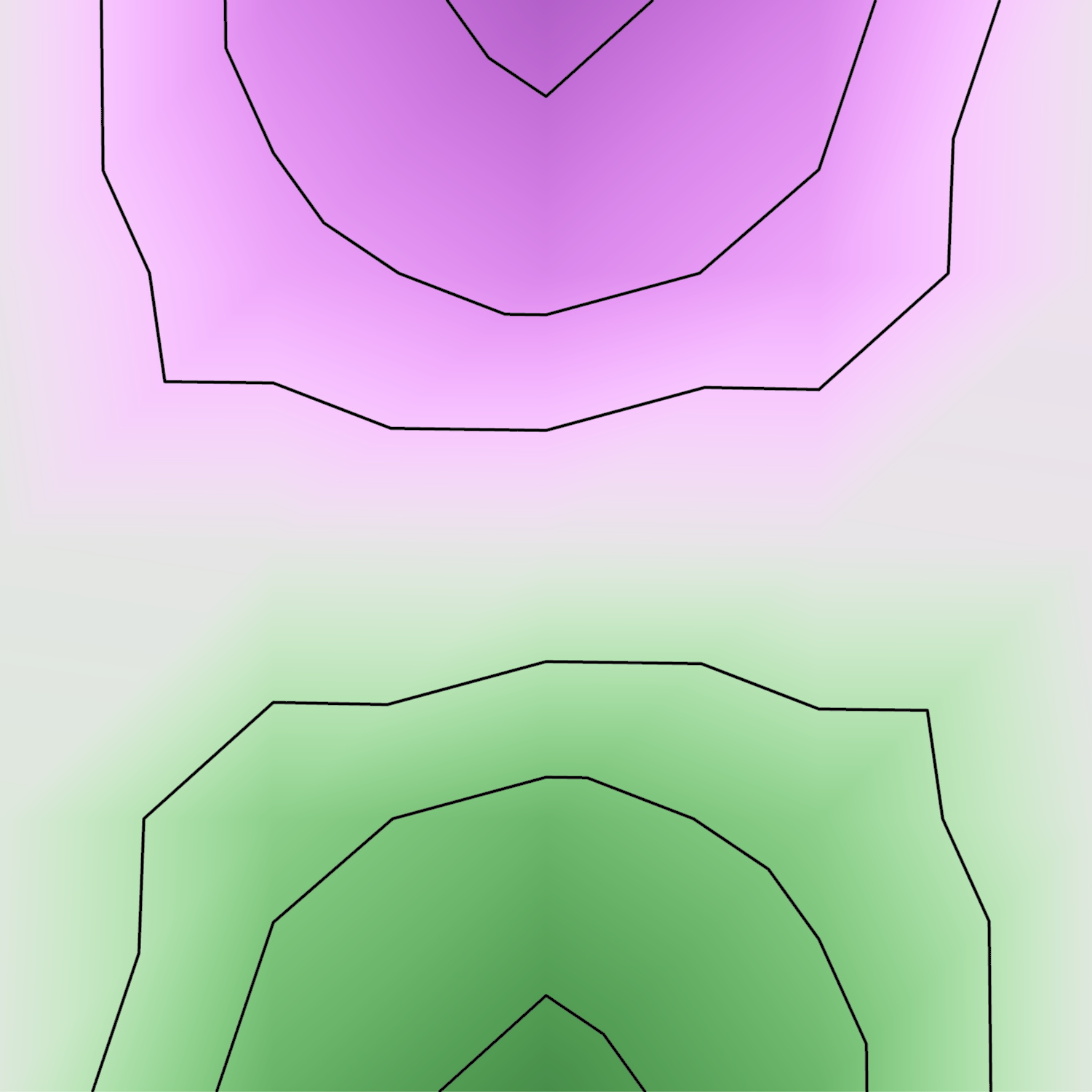}}
\caption{Plot of the $x$-component of the solution velocity
of \cite[(18.19)]{lrsBIBih} using BDM(1):
(a) discontinuous result, (b) after filtering using SISG.}
\label{fig:mixed}
\end{figure}

One approach to approximating flow in porous media with discontinuous 
physical properties is to use
mixed methods \cite{lrsBIBih} and discontinuous finite elements.
This approach is called discontinuous Galerkin (DG).
Thus visualization of the primary velocity variable is a challenge due to
its discontinuity.
DG methods are also widely applied to hyperbolic PDEs \cite{cockburn2003enhanced}. In Figure \ref{fig:mixed}, we depict the solution of the problem in \cite[(18.19)]{lrsBIBih} using BDM elements of order 1.

\begin{figure}
\centerline{\includegraphics[width=6in]{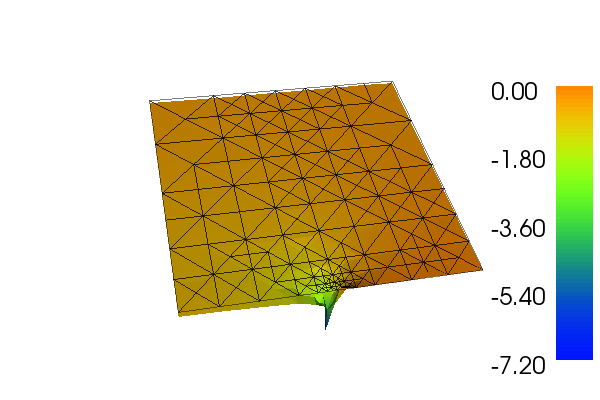}}
\caption{Derivative of singular function: SISG in {\tt DOLFIN} \cite{logg2010dolfin}.
Solution of the problem in Section \ref{sec:hirefme} with tolerance 0.002.}
\label{fig:adaptikslip}
\end{figure}

\subsection{Highly refined meshes}
\label{sec:hirefme}

Highly refined meshes are used to resolve solution singularities in many
contexts.
For example, they can be required to resolve singularities in data
\cite[Section 6.2]{lrsBIBih}.
Even when data is smooth, singularities can arise due to domain geometry or 
changes in boundary condition types \cite[Section 6.1]{lrsBIBih}.
Such singularities are common, and physical quantities are often associated
with derivatives of solutions to such problems.
For example, the function $g$ defined in polar coordinates by
$$
g(r,\theta)=r^{1/2}\sin(\half\theta)
$$
arises naturally in this context.
Consider the boundary value problem to find a function $u$ satisfying
the partial differential equation
$$
-\Delta u = 1
$$
in the domain $\Omega:=[-\half,\half]\times[0,1]$, 
together with the boundary conditions
$$
u=0 \;\hbox{on}\; \set{(x,0)}{0\leq x\leq \half}, \quad
\frac{\partial u}{\partial n}=0 \;\hbox{on}\;\set{(x,0)}{-\half\leq x\leq 0},
$$
with $u=g-\half y^2$ on the remainder of $\partial\Omega$.
Note that $g$ is a smooth function on this part of the boundary.
Then the exact solution is $u=g-\half y^2$.

This problem has a singularity at $(0,0)$;
the derivative $u_x$ is infinite there.
The solution $u$ of this problem is shown in \cite[Figure 16.1]{lrsBIBih}, 
computed using automatic, goal-oriented refinement and 
piecewise-linear approximation.
Here we focus instead on $u_x$.
In this case, the derivative of the piecewise linear approximation
is a piecewise constant, and so the SISG approach seems warranted.
Thus we project $u_x$ onto continuous piecewise linear functions on the same
mesh, and the result is depicted in Figure \ref{fig:adaptikslip}.
This approach is actually the default approach taken in {\tt DOLFIN}
\cite{logg2010dolfin}.
We see that SISG is quite effective in representing a very singular, 
discontinuous computation in a comprehensible way.

To make this result more quantitative, we compare with the exact solution
derivative $u_x$. 
Note that we can write
$$
g(x,y)=(x^2+y^2)^{1/4} \sin\big(\half\hbox{atan}(y/x)\big),
$$
so that
$$
g_x(x,y)=\half(x^2+y^2)^{-3/4} \big(x\sin(\half\hbox{atan}(y/x))
- y\cos(\half\hbox{atan}(y/x))\big).
$$
To avoid the singularity at the origin in the computations, we introduce $\epsilon>0$
and replace $g_x$ by
$$
g^\epsilon_x(x,y)=\half(\epsilon+x^2+y^2)^{-3/4} \big(x\sin(\half\hbox{atan}(y/x))
- y\cos(\half\hbox{atan}(y/x))\big).
$$
This introduces a small error that we minimize with respect to $\epsilon$.
More precisely, we compute
\begin{equation}\label{eqn:whatecomput}
\norm{\Pi (u_{h,x})-g^\epsilon_x}_{L^2(\Omega)}.
\end{equation}
The computational results are given in Table \ref{tabl:exactickslip}.
The goal of the adaptivity was to minimize
\begin{equation}\label{eqn:stdbndsmput}
\norm{\nabla (u_{h}-u)}_{L^2(\Omega)},
\end{equation}
and the data in Table \ref{tabl:exactickslip} indicates that this was successful.
\color{red}
We refer to \eqref{eqn:stdbndsmput} as the ``grad error''
and it is presented in the 4-th column in Table \ref{tabl:exactickslip}.
\color{black}
The initial mesh consisted of eight $45^\circ$ right-triangles, with 9 vertices,
in all of the computations.
For tolerances greater than 0.03, no adaptation occurs.

\color{red}
The values in \eqref{eqn:stdbndsmput} 
\color{black}
were computed via the {\tt DOLFIN} function {\tt errornorm} 
without any regularization of the exact solution $u(x,y)=g(x,y)-\half y^2$.
This quantity is a non-SISG error, in that the gradients are treated as
discontinuous piecewise-defined functions.
By contrast the SISG error \eqref{eqn:whatecomput} corresponds to part
of the quantity in \eqref{eqn:stdbndsmput} (the $x$-derivative), but this error
is a measure of the accuracy of the SISG projection.

\begin{figure}
\centerline{\includegraphics[width=5.5in]{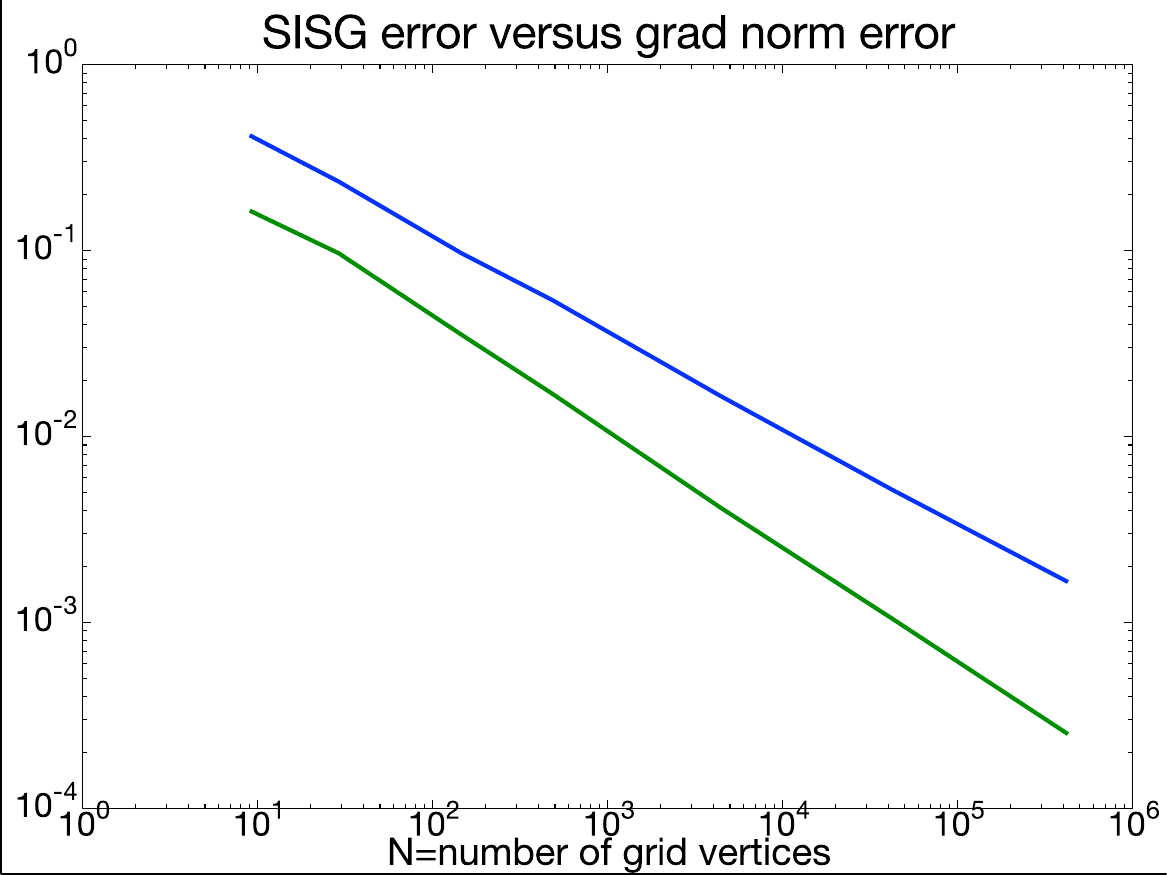}}
\caption{
\color{red}
Comparison of the error \eqref{eqn:whatecomput} for the SISG technique and
the grad error \eqref{eqn:stdbndsmput} as a function of
N, the number of vertices in the mesh, including boundary
vertices, after adaption.}
\label{fig:revslipa}
\end{figure}

Both \eqref{eqn:stdbndsmput} and \eqref{eqn:whatecomput} decay approximately like 
$CN^{-1/2}$, where $N$ is the number of vertices in the mesh (including boundary
vertices) after adaptation, 
\color{red}
as indicated in Figure \ref{fig:revslipa} and Table \ref{tabl:exactickslip}.
Indeed, the rate of decay for the SISG error appears to be slightly faster.
\color{black}
For a smooth $u$ and a regular mesh of size $h$, we would expect the quantities
\eqref{eqn:stdbndsmput} and \eqref{eqn:whatecomput} to be $\order{h}$,
and $N=\order{h^{-2}}$ in this case.
So the observed convergence $CN^{-1/2}$ in Table \ref{tabl:exactickslip}
is best possible.
Note that the tolerance value is associated with the square of the
quantity \eqref{eqn:stdbndsmput}.

\begin{table}[htbp]
\centering
\begin{tabular}{|r|r||c|c||c|c||c|c|}
\hline
N & tolerance & SISG error \eqref{eqn:whatecomput} & rate & grad error  & rate& $\epsilon$ & time\\
\hline
9      &  0.03   & 1.64e-01 & & 4.17e-01& & 1.00e-03 & 0.016\\
29     &  0.01   & 9.68e-02 &0.45& 2.35e-01&0.49& 1.00e-03 & 0.017 \\
146    & 0.002   & 3.53e-02 &0.62& 9.68e-02&0.55& 1.00e-05 & 0.019 \\
482    &  0.001  & 1.70e-02 &0.61& 5.42e-02&0.49& 1.00e-06 & 0.019 \\
4477   & 0.0001  & 4.11e-03 &0.64& 1.64e-02&0.54& 1.00e-08 & 0.038 \\
41481  & 0.00001 & 1.06e-03 &0.61& 5.22e-03&0.51& 1.00e-11 & 0.223\\
427169 & 0.000001 & 2.51e-04&0.62& 1.65e-03&0.49& 1.00e-13 & 2.089\\
\hline
\end{tabular}
\caption{Errors \eqref{eqn:whatecomput} for the SISG technique,
\color{red}
based on adaptive, piecewise-linear, finite-element computation.
\color{black}
The initial mesh size was 2 for all of the computations.
N denotes the number of vertices in the mesh (including boundary
vertices) after adaptation,
the second column indicates the tolerance used for adaptivity, 
the third column denotes the error quantity defined
in \eqref{eqn:whatecomput}, 
\color{red}
``grad error'' denotes \eqref{eqn:stdbndsmput}, and
\color{black}
the last column gives the value of $\epsilon$ used in \eqref{eqn:whatecomput}.
\color{red}
The ``rate'' columns give the rate $\rho$ of decrease for the errors modeled as $N^{-\rho}$.
The last column ``time'' gives the time in seconds to compute the SISG projection using {\tt dolfin} on a single core of a MacBook Pro with a 2.3 GHz Intel Core i7.
\color{black}
}
\label{tabl:exactickslip}
\end{table}

\section{Methods}

Except for Figure \ref{fig:adaptikslip}, the images presented here were 
produced using Firedrake \cite{rathgeber2017firedrake}
and Diderot \cite{Chiw-Diderot-PLDI-2012,Kindlmann-Diderot-VIS-2015}.
Figure \ref{fig:adaptikslip} was done with DOLFIN \cite{logg2010dolfin}.

\section{Conclusions}

The smoothing technique of Savitsky and Golay can be extended 
to finite element methods in a useful way.
This allows accurate presentation of derivatives of piecewise-defined 
functions, even for highly refined meshes.
It also allows discontinuous approximations, such as Discontinuous Galerkin,
to be visualized in an effective way.



\clearpage
\bibliographystyle{siamplain}
\bibliography{siac}
\end{document}